\documentclass{amsart}

\usepackage{amssymb,amsmath,amsthm,amscd}%--pas necessaire avec amsart

\usepackage[active]{srcltx} % ajout Kristian
\usepackage[frenchb]{babel}
\usepackage{pstricks}
\usepackage[T1]{fontenc}
\usepackage{mathrsfs}

\newcommand{\C}{\mathbb{C}}%------------------------Complex numbers
\newtheorem{theo}{Th\'eor\`eme}[section]%
\newtheorem{prop}[theo]{Proposition}%
\newtheorem{rem}[theo]{Remarque}%
\newtheorem{lem}[theo]{Lemme}%
\newtheorem{defn}[theo]{D\'efinition}%

\begin{document}

\title[Bohr's phenomenon on a regular condensator in the complex plane]{Bohr's phenomenon on a regular condensator in the complex plane}
\date{\today}

\author[P. LassËre]{LassËre Patrice}
\email{lassere@math.ups-tlse.fr}
\address{LassËre Patrice : Institut de MathÈmatiques, , UMR CNRS 5580,
Universit\'e Paul Sabatier, 
118 route de Narbonne, 31062 TOULOUSE, FRANCE}
%delete this block if only one author, copy if more than two authors
\author[E. Mazzilli]{Mazzilli Emmanuel}
\email{Emmanuel.Mazzilli@math.univ-lille1.fr}
\address{UniversitÈ Lille 1, 59655 Cedex, VILLENEUVE D'ASCQ,  FRANCE.}


% fill out if necessary or keep empty, acknowledgements go before bibliography
%\thanks{Research supported by XYZ}
%\dedicatory{}

\begin{abstract}
  We prove the following generalisation of Bohr's theorem : let $K\subset\mathbb C$ a continuum, $(F_{K,n})_{n\geq 0}$ its Faber polynomials, $\Omega_R$ the level sets of the Green function of $\bar{\mathbb C}\setminus K$ with singularity at infinity,  then there exists $R_0$ such that for any $f=\sum_n a_n F_{K,n}\in\mathscr O(\Omega_{R_0})$ : $f( \Omega_{R_0})\subset D(0,1)$ implies  $\sum_n\left\vert a_n \right\vert\cdot\Vert F_{K,n}\Vert_K<1$.
\end{abstract}

%%%%%%%%%%%%%%%%%%%%%%%%%%%%%%%%%%%%%%%%%%%%%%%%%%%%%%%%%%%%%%%%
%opening

\maketitle

\section{Introduction}

\bigskip

The well-known  Bohr's theorem  \cite{bohr}  states that for any function   $f(z)=\sum_{n\geq 0}\,a_n z^n$ holomorphic on the unit disc $\mathbb D$ : 
$$\left(\  \left \vert\sum_{n\geq 0}\, a_n z^n\right\vert<1,\ \forall\,z\in\mathbb D\ \right)\ \implies\ 
\left( \sum_{n\geq 0}\, \left\vert a_n z^n\right\vert<1,\ \forall\,z\in D(0,1/3)\ \right)$$
and the  constant $1/3$ is optimal. 

\bigskip
Our goal in  this work is to study   Bohr's theorem  in the following context. 
Let $K\subset\mathbb C$ be a compact in the complex plane. What are the open sets $\Omega$  containing $K$ such that the space $\mathscr O(\Omega)$ admits a topological basis\footnote{For all $f\in\mathscr O(\Omega)$ there exists an unique sequence $(a_n)_n$ of complex numbers such that  $f=\sum_{n\geq 0} a_n\varphi_n$ for the usual compact convergence topology of $\mathscr O(\Omega)$.}   $(\varphi_n)_n$ which  verifies, for every holomorphic function $f=\sum_{n\geq 0} a_n\varphi_n\in\mathscr O(\Omega)$ :
$$\left(\  \left \vert\sum_{n\geq 0}\, a_n \varphi_n(z)\right\vert<1,\ \forall\,z\in\Omega\ \right)\ \implies\ 
\left( \sum_{n\geq 0}\, \left\vert a_n \right\vert\cdot\Vert\varphi_n\Vert_K<1\ \right)\ ?$$
In this case we say that the family $(K, \Omega, (\varphi_n)_{n\geq 0})$ satisfies  \textbf{Bohr's property} or that \textbf{Bohr's phenomenon} is observed.

\bigskip
\noindent \textbf{Some examples : } $\bullet$ The family $(\overline{D(0,1/3)}, D(0,1), (z^n)_{n\geq 0})$ satisfies  Bohr's phenomenon (this is Bohr's classic theorem).

\noindent $\bullet$    Note that the family $(\overline{D(0,1/3)}, D(0,1), ((3z)^n)_{n\geq 0})$ also satisfies  Bohr's phenomenon. This example will play a special role in the following, since  $((3z)^n)_{n\geq 0}$ is the Faber polynomial basis  associated with the compact $\overline{D(0,1/3)}$.

\noindent $\bullet$ On the other hand, the family $(\overline{D(0,2/3)}, D(0,1), (z^n)_{n\geq 0})$ does not satisfy   Bohr's phenomenon (due to optimality of the constant  $1/3$ in Bohr's theorem).

\bigskip
As a starting point, for a given compact $K$  we must choose  a ``good'' open neighborhood $\Omega$, that  admits for  $\mathscr O(\Omega)$  a ``nice''  basis  $(\varphi_n)_n$. ``Nice''  here means that there are good local estimates for  $\varphi_n$ on $\Omega$ but not only, since, unlike for other well-known theorems for power series on the disc  \cite{lasserenguyen},   Bohr's theorem cannot be extended to all basis. For example, as pointed out by Aizenberg  \cite{AAD}, it is necessary that one of the elements of the basis be a constant function. 

\bigskip
We want to focus on the following situation :

\bigskip
\begin{defn} Let $K$ be a compact in $\mathbb C$ including at least two points, $K$ is a continuum if  $\overline{\mathbb C}\setminus K$ is simply connected. 
\end{defn}

\bigskip
When  $K$ is a continuum it can be associated with the sequence $(F_{K,n})_n$ of its Faber polynomials. In more detail, let $\Phi\ :\ \overline{\mathbb C}\setminus K\to \overline{\mathbb C}\setminus{\overline {\mathbb D}}$ be the unique conformal mapping that verifies
$$\Phi(\infty)=\infty,\quad \Phi'(\infty)=\gamma>0.$$
Therefore $\Phi$ admits a Laurent development close to the infinity point under the form:
$$\Phi(z)=\gamma z+\gamma_0+\dfrac{\gamma_1}{z}+\dots+\dfrac{\gamma_k}{z^k}+\dots$$
and then for $n\in\mathbb N$ :
$$\begin{aligned}\Phi^n(z)&=\left( \gamma z+\gamma_0+\dfrac{\gamma_1}{z}+\dots+\dfrac{\gamma_k}{z^k}+\dots \right)^n\\
&=\underbrace{\gamma^n z^n+a_{n-1}^{(n)}z^{n-1}+\dots+a_{1}^{(n)}z+a_{0}^{(n)}}_{ F_{K,n}(z)}+\underbrace{  \dfrac{b_{1}^{(n)}}{z}+\dfrac{b_{2}^{(n)}}{z^2}+\dots+\dfrac{b_{k}^{(n)}}{z^k}+\dots}_{ E_{K,n}(z)} 
\end{aligned}$$
$F_{K,n}$ is   the polynomial part of the Laurent expansion at infinity of $\Phi^n$. It is a common basis for the spaces  $\mathscr O(K),\ \mathscr O(\Omega_R),  (R>1)$ where\footnote{\samepage $\Omega_R$ is also the level set of the Siciak-Zaharjuta extremal function $\Phi_K(z):=\sup\{\vert p(z)\vert^{1/{\text{deg}}(p)}\}$ where the supremum is taken over all complex  polynomials $p$ such that $\Vert p\Vert_K\leq 1$. $\Phi_K$ is also related to the classical Green function for $\bar{\mathbb C}\setminus K$   with pole at infinity $g_K\ :\ \mathbb C\setminus K\to ]0,+\infty[$ by the equality $\log\Phi_K=g_K$ on $\mathbb  C\setminus K$ . Recall that  $g_K$ is the unique harmonic positive function on $\mathbb C\setminus K$ such that $\lim_{z\to\infty} \left( g_K(z)-\log\vert z\vert\right) $ exists and is finite and $\lim_{z\to w}g_K(z)=0,\ \forall\,w\in\partial(\mathbb C\setminus K)$. } $\Omega_R:=\{ z\in\mathbb C\ :\ \vert \Phi(z)\vert<R\}\cup K$. This polynomial basis exhibits remarkable properties (the relevant reference is the work by P.K.Suetin \cite{suetin}) similar to the Taylor basis $(z^n)_n$ on discs $D(0,R)$. In particular, the level sets $\Omega_R$ are the convergence domains of the series $\sum_{n\geq 0} a_n F_{K,n}$ and for any compact $L\subset \overline{\mathbb C}\setminus K$ we have
$$\lim_{n\to\infty} \Vert F_{K,n}\Vert_L^{1/n}= \Vert \Phi\Vert_L. $$
This formula is the one variable version of a more general formula (see \cite{nguyen}).

\bigskip
In this work, we show (Theorem 3.1) that for every continuum $K$ there exists an $R_0>1$ such that for any $R\geq R_0$ the family 
$(K,\Omega_R, (F_{K,n})_{n\geq 0})$ verifies  Bohr's property.
 We start  by studying the cases of an  elliptic condensator (i.e. $K=[-1,1]$) which had been considered in a different form by Kaptanoglu  and Sadik   in an interesting study  \cite{kap} which motivated this article (see  remark 2.4).

 \bigskip
\noindent \textbf{Acknowledgement.} Finally, we thank the Anonymous Referee for useful suggestions  
improving significantly the paper.

\bigskip

\section{An example : the  ``elliptic'' condensator  $K=[-1,1]$}

\bigskip
 Let us examine in this section the particular case where  $K:=[-1,1 ]$. This is a ``fundamental'' example because this is one of the very few case (see \cite{suetin}, \cite{he} for circular lunes) where the explicit form of  the conformal map $\Phi \ :\ \Omega:=\overline{\C}\setminus K\to \{\vert w\vert >1\}$ allows us to obtain a more precise estimation of the Faber polynomials of $K$  (see \cite{suetin}).

\noindent  Here, $\Phi^{-1}(w)={1\over 2}(w+{w}^{-1})$ is the Zhukovskii function, the  Faber polynomials $(F_{K,n})_n$ form a common basis for the spaces  $\mathscr O(\Omega_R)$, $(R>1)$ where the boundary $\partial\Omega_R=\Phi^ {-1}(\{\vert w \vert=R\}$ of the level set  $\Omega_R$ is given by the equation :
$$2z=R e^{i\theta}+R^{-1} e^{-i\theta}.$$
Theses are ellipses with foci $1$ et $-1$ and eccentricity $\varepsilon= \frac{2R}{1+R^2}$.
We observe that  the polynomials $F_{K,n}$ enjoy in the target coordinates ``$w$'' a much more convenient form for computation than in the source coordinates ``$z$''.
Indeed  $\Phi$ presents a simple pole at infinity which implies that
$\Phi^n+{1/ \Phi^n}$ et $\Phi^n$ have the same principal part. We observe also that
$\Phi(z)=z+\sqrt{z^2-1},$
which implies
${1/\Phi(z)}=z-\sqrt{z^2-1}$. From these last identities we can deduce\footnote{We can also deduce (see \cite{suetin}, pp. 36-37) that if $K=[-1,1]$, then the Faber polynomials are the Tchebyshev polynomials of the first kind (up to a constant $2$ if $n\geq 1$)  : $F_{K,0}(z)=T_0(z),\ F_{K,n}(z)=2T_n(z),\ (n\geq 1)$ where $T_n(x)=\cos(n{ \text{arccos}} x)$.} that
 ${1/\Phi^n}+\Phi^n$ extends as a polynomial  on $\Bbb{C}$. This is $F_{K,n}$  and if we write  $F_{K,n}$ in the  target coordinates ``$w$'', we get :
$$F_{K,n}(w)=w^n+ w^{-n}.$$
This important equality will allow us to write any function $f(z)=\sum_n a_n F_{K,n}(z), z\in \Omega_R$, holomorphic on $ \Omega_R$ under the form
$$f(z)=f(\Phi^{-1}(w))=\sum_n a_n F_{K,n}((\Phi^{-1}(w))=\sum_n a_n \left( w^n+w^{-n}\right),\quad 1<\vert w\vert <R,$$
and we shall often use this device from now on.

 \bigskip
 Now let us look at Bohr's phenomenon for the elliptic condensator 
 $(K:=[-1,1], \Omega_R, (F_{K,n})_{n\geq 0})$ given that $R>R_0$
is large enough. 
 Then next proposition is,  in our particular case, 
the equivalent version of Caratheodory's inequality.

\begin{prop} Let
$f(w)=a_0+\sum_{1}^{\infty}a_n(w^n+ w^{-n})\in \mathscr O(\{1<\vert w\vert <R\})$. Suppose that
$\texttt{re}({f})>0$,  then : 
$$\vert a_n\vert\leq {2\texttt{re}(a_0)\over
R^n- R^{-n}},\quad\forall\,n>0. $$
\end{prop}

\bigskip
\noindent\textbf{Proof : } Let  $1<r<R$, then  for all $n>0$  we have
 $$\begin{aligned}&a_n r^{-n}&=&{1\over
2\pi}\int_{0}^{2\pi}e^{in\theta}f(r e^{i\theta})d\theta,\\
&\overline a_nr^n&=&{1\over
2\pi}\int_{0}^{2\pi}e^{in\theta}\bar f(r e^{i\theta})d\theta.
\end{aligned}$$
 which easily gives
(remember that $\texttt{re}(f)>0$) :
$$\vert a_n\vert \cdot\left(r^n-r^{-n}\right)\leq\left\vert {a_nr^{-n}}+\bar a_nr^n\right\vert\leq {1\over
\pi}\int_{0}^{2\pi}\texttt{re}(f(r e^{i\theta}))d\theta=2\texttt{re}({a_0}),$$
to get the expected result, (just let $r$ tend to $R$).\hfill$\blacksquare$
\medskip

\begin{lem} Let $f=a_0+\sum_{n=1}^{\infty}a_n(w^n+w^{-n}) \in \mathscr O(\{1<\vert w\vert <R\})$. Suppose that $\vert
f\vert <1$ and $a_0>0$, then\footnote{ Note that $\vert f\vert<1$ implies $a_0<1$.} we have :
$$\vert a_n\vert\leq {2(1-a_0)\over R^n-R^{-n}}.$$
\end{lem} 

\bigskip
\noindent\textbf{Proof : } This is classical : let $g=1-f$, then
$\texttt{re}({g})>0$ on $\{1<\vert w\vert <R\}$ and by prop. 2.1 :
$$\vert a_n\vert\leq {2(1-a_0)\over R^n-R^{-n}}.$$
\hfill$\blacksquare$

\begin{prop} For all $R\geq R_0=5.1284...$ the family  $(K:=[-1,1], \Omega_R, (F_{K,n})_{n\geq 0})$ satisfies Bohr's phenomenon ($ \Omega_{R_0}$ is the ellipse with eccentricity $\varepsilon_0=0.3757...$). \end{prop}

\bigskip
\noindent\textbf{Proof : } Let  $f=a_0+\sum_{1}^{\infty}a_nF_{K,n}\in\mathscr O(\Omega_R)$ and suppose that $\vert f\vert <1$ on
$\Omega_{R}$. In the variables ``$w$'' :
$f(w)=a_0+\sum_{1}^{\infty}  a_n (w^n+ w^{-n})$ on $\{1<\vert
w\vert <R\}$ and up to a  rotation (changing nothing by symmetry), we can suppose that
$a_0\geq 0$. Then by lemma 2.2 :
$$\begin{aligned}a_0+\sum\vert a_n\vert\cdot\Vert
F_{K,n}\Vert_{K}&\leq
a_0+2(1-a_0)\sum_{n=1}^{\infty}{r^n+{r^{-n}}\over
R^n-{R^{-n}}},\quad (1<r<R)\\
&\leq a_0+(1-a_0)\sum_{n=1}^{\infty}{4R^n\over R^{2n}-1}.
\end{aligned}$$
%On veut montrer que pour $R$ assez grand, il existe $1<r<R$ tel que:
%$$2(1-a_0)\sum_{1}^{\infty}{r^n+{r^{-n}}\over
%R^n-{R^{-n}}}<1-a_0,$$ 
%ceci est possible si et
%seulement si $\inf_{[1,R]}\phi_{R}<1$, o\`u $\phi_{R}$
%est la fonction d\'efinie par:
%$$\phi_{R}(r):=2\sum_{1}^{\infty}{r^n+{r^{-n}}\over
%R^n-{R^{-n}}}.$$ Cette fonction \'etant
%croissante sur $[1,R]$, il faut donc que 
This gives
$$a_0+\sum_{n=1}^{\infty}\vert a_n\vert\cdot\Vert
F_{K,n}\Vert_{K}<1$$
if  
$$\varphi(R):=\sum_{1}^{\infty}{4R^n\over R^{2n}-1}<1.$$ 
But $\varphi$ strictly decreases on $]1,\infty[$, $\lim_{1_+}\varphi(R)=+\infty$,  $\lim_{+\infty}\varphi(R)=0$ therefore, there exists a unique 
$R_0>1$  such that $\varphi(R)-1=0$ on
$]1,\infty[$ ;  Mathematica  gives $R_0=5.1284...$ corresponding to an eccentricity of $\varepsilon_0=0.3757...$ ;   $(K:=[-1,1], \Omega_R, (F_{K,n})_{n\geq 0})$ satisfies Bohr's phenomenon for all  $R\geq R_0$.\hfill$\blacksquare$

\bigskip
\begin{rem}    Using theorem 7 in \cite{kap}, we can deduce  a weaker version of proposition 2.3 with $R_0=5.1573...$ and $\varepsilon_0=0.3738..$, so,  proposition 2.3 is a slighty stronger version of  theorem 7 in \cite{kap}. In another work \cite{lasseremanu} we calculate exactly  the infimum of $R_0$ satisfying proposition  2.3 i.e. what we call the Bohr's radius of $K=[-1,1]$ in Theorem 3.1.
\end{rem}

\medskip

\section{ Bohr's phenomenon on an arbitrary Green condensator}

\bigskip
\subsection{Estimations of Faber polynomials on a Green condensator } In this paragraph, we recall classical inequalities (see \cite{suetin}) on Faber polynomials of $K$ that we will use in paragraph 3.2. 

\medskip
Let $K\subset\mathbb C$ be a continuum, $(F_{K,n})_{n\geq 0}$ its Faber polynomials. Recall that $\Phi^n(z)=F_{K,n}(z)+E_{K,n}(z)$ where $E_{K,n}$ is the meromorphic part in the Laurent developement of $\Phi^n$ in a neighborhood of infinity. If  $\Omega_r$, $(r>1)$ is the level set $\{ z\in\mathbb C\ :\ \vert\Phi(z)\vert<r\}$ then  we have the following integral formulas for Faber polynomials (see Suetin, \cite{suetin}, pp 42) : 
$$\forall\,z\in\Omega_r\ :\quad F_{K,n}(z)=\int_{\partial\Omega_r}{\Phi^n(t)\over t-z}dt,\ \ \leqno{(1)}$$
$$\forall\,z\in\mathbb C\setminus \overline{\Omega_r}\ :\quad E_{K,n}(z)=\int_{\partial\Omega_r}{\Phi^n(t)\over t-z}dt,\ \ \leqno{(2)}$$
Formula (2) leads to the following estimations for all $1<r<R$ :
$$\forall\,z\in\mathbb C\setminus {\Omega_R}\ :\quad\vert E_{K,n}(z)\vert\leq \int_{\partial\Omega_r}\left\vert{\Phi^n(t)\over  t-z}\right\vert\cdot \vert
dt\vert\leq {r^n\texttt{lg}(\partial\Omega_r)\over \texttt{dist}(z,\partial\Omega_r)},\ \leqno{(3)}$$
($\texttt{lg}(\partial\Omega_r)$ is the euclidian length $\partial\Omega_r$, 
$\texttt{dist}(z,\partial\Omega_r)$ is the euclidian distance from $z$ to $\partial\Omega_r$)  and
$$\forall\,z\in\partial\Omega_R\ :\ \quad \vert F_{K,n}(z)\vert\leq R^n\left(1+{r^n\over R^n}\cdot {\texttt{lg}(\partial\Omega_r)\over \texttt{dist}(z,\partial\Omega_r)}\right)
,\ \ \leqno{(4)}$$
 for all $1<r<R$. Then if $R$ is large enough, precisely if  $${r^n\over
R^n}\cdot{\texttt{lg}(\partial\Omega_r)\over \texttt{dist}(z,\partial\Omega_r)}<1$$ 
then for all $n>0$, we have :
$$\forall\,z\in\partial\Omega_R\ :\ \quad\vert F_{K,n}(z)\vert\geq R^n\left(1-{r^n\over R^n}\cdot{\texttt{lg}(\partial\Omega_r)\over
\texttt{dist}(z,\partial\Omega_r)}\right)>0.$$
With formula (1), we deduce the estimation,  for all $ r>1$ and $z\in
K$ :
$$\vert F_{K,n}(z)\vert\leq \int_{\partial\Omega_{r}}\left\vert{\Phi^n(t) \over  t-z}\right\vert\cdot\vert
dt\vert\leq r^n{\texttt{lg}(\partial\Omega_{r})\over \texttt{dist}(z,\partial\Omega_{r})}.\ \leqno{(5)}$$
 
 If moreover the compact $K$ is a domain defined by a real analytic Jordan curve, then Caratheodory's theorem  ensures that $\Phi$ extends as a biholomorphism on a neighborhood  of $\partial K$, say up to  $\partial\Omega_{r_0}$, where $r_0<1$. From this, we get for all $r_0<R$ :
$$\forall\,z\in\mathbb C\setminus\Omega_R\ :\ \quad\vert E_{K,n}(z)\vert\leq \int_{\partial\Omega_{r_0}}\left\vert{\Phi^n(t)\over  t-z}\right\vert\cdot\vert
dt\vert\leq {r_{0}^n\texttt{lg}(\partial\Omega_{r_0})\over \texttt{dist}(z,\partial\Omega_{r_0})},$$ 
and so the estimations
$$\begin{aligned} R^n\left(1-{r_{0}^n\over R^n}\cdot{\texttt{lg}(\partial\Omega_{r_0})\over \texttt{dist}(z,\partial\Omega_{r_0})}\right)&\leq\vert F_{K,n}(z)\vert\\
&\leq R^n\left(1+{r_{0}^n\over R^n}\cdot
{\texttt{lg}(\partial\Omega_{r_0})\over \texttt{dist}(z,\partial\Omega_{r_0})}\right),
 \end{aligned}$$
for all $z\in\mathbb C\setminus\Omega_R,\ r_0<R$.
\bigskip

\bigskip
\subsection{Bohr's phenomenon on a Green condensator. } In this  paragraph we extend  proposition 2.3  for all continuum  $K$ in the complex plane, precisely :

\bigskip
\begin{theo}
For all continuum $K\subset\mathbb C$, there exists a constant $R_K>1$ such that for all $R>R_K$ the family $(K,\Omega_R,(F_{K, n})_{n\geq 0},)$ satisfies  Bohr's phenomenon and the infimum $R_0$ of such $R$ will be called    the \textbf{Bohr's radius} of $K$.
\end{theo}

%\bigskip
%Then we'll define the Bohr's radius $R_0$ of a compact $K$ as the infimum of the $R>1$ such that the family $(K,\Omega_R,(F_{K, n})_n,)$ vÈrifie le phÈnomËne de Bohr. 

For example the Bohr radius for a disc  $K=D(a,r)$ is $3$ due to Bohr's classic theorem, and in \cite{lasseremanu} we compute the exact value  of $R_0$ when $K=[-1,1]$.

\bigskip
Before proving theorem 3.1, some intermediate results are necessary. Let $K$ be a  continuum, $(F_{K,n})_{n\geq 0}$ its sequence of Faber polynomials and $z_0\in \partial K$. Consider the family $(\varphi_{n\geq 0})_{n\geq 0}$ where   $\varphi_0\equiv 1$ and
$\varphi_n=F_{K,n}-F_{K,n}(z_0)\  (n\geq 1)$. It is clear that  $(\varphi_n)_{n\geq 0}$ is again a basis of the spaces $\mathscr O(\Omega_R)$ for all $R>1$ and we have

\begin{theo} The family $( K, \Omega_R, (\varphi_n)_{n\geq 0})$ enjoys  Bohr's property for  $R$ large enough. That is to say, there exists $R>1$ such that   all holomorphic function $f=\sum_n a_n\varphi_n\in\mathscr O(\Omega_R) $ with values in  $\mathbb D$ satisfy 
$$\sum_{n\geq 0} \vert a_n\vert\cdot\Vert \varphi_n\Vert_K=\vert f(z_0)\vert+\sum_{n\geq 1} \vert a_n\vert\cdot\Vert \varphi_n\Vert_K <1.$$
\end{theo}

\bigskip
\noindent \textbf{Proof : } Let $R_0>1$. We can suppose without loss of generality  that $z_0=0$. Because  $\varphi_n(0)=0$ for all $n>0$ we can apply 
theorem 3.3 in \cite{AAD} on the open set $\Omega_{R_0}$. This implies that there exists  
$D(0,\rho_0)$ where $\rho_0$ is small enough and a compact $K_1\subset\Omega_{R_0}$ such that :
$$\vert f(0)\vert
+\sum_{n\geq 1}\vert a_n\vert\cdot \Vert
\varphi_n\Vert_{D(0,\rho_0)}\leq \Vert f\Vert_{K_1},$$
for any function $f=\sum_n\, a_n\varphi_n\in\mathscr O(\Omega_{R_0})$.
Now choose  $\rho_1>0$ such that $K_1\subset D(0,\rho_1)$. We have :
$$\vert f(0)\vert
+\sum_{n\geq 1}\vert a_n\vert\cdot \Vert
\varphi_n\Vert_{D(0,\rho_0)}\leq \Vert f\Vert_{D(0,\rho_1)},\leqno{(6)}$$ 
for all $f=\sum_n a_n\varphi_n\in\mathscr O(\Omega_R)$ where $R$ is choosen large enough so that $D(0,\rho_1)\subset \Omega_R$.

\medskip
Let $f\in \mathscr O(\Omega_R)$ such that  $\Vert f\Vert_{\Omega_R}\leq 1$; the invariant form of  Schwarz's lemma (\cite{goluzin}, chapter 8) gives the following estimation on any disc 
$D(0,\rho)\subset \Omega_R$ ($\rho\geq \rho_1$) :
$$\Vert f\Vert_{D(0,\rho_1)}\leq {\rho_1 \rho^{-1}+\vert
f(0)\vert\over 1+\vert f(0)\vert \rho_1 \rho^{-1}}.\leqno{(7)}$$
We want for $f=f(0)+\sum_{n\geq 1}\, a_n\varphi_n\in\mathscr O(\Omega_R)$
to dominate the quantity :
$\vert f(0)\vert+\sum_{n\geq 1}\, \vert a_n\vert\cdot\Vert\varphi_n\Vert_{K}$; write
$$\sum_{n\geq 1}\vert
a_n\vert\cdot\Vert\varphi_n\Vert_{K}=\sum_{n\geq 1}\vert
a_n\vert\cdot\Vert\varphi_n\Vert_{D(0,\rho_0)}\times
\dfrac{\Vert \varphi_n \Vert_{K}}{ \Vert \varphi_n \Vert_{D(0,\rho_0)}}.\ \
\leqno{(8)}$$
Let $L$ be a disc contained in $D(0,\rho_0)\setminus K$ then
$$\lim_{n\to\infty} \Vert\varphi_n\Vert_L^{1/n} =R^{\alpha_L}$$ 
where\footnote{$\omega$ is the extremal  function associated for the pair $(K, \Omega_R)$.} $\alpha_L:=\max_{z\in L} \omega(z, K,\Omega_R)$, this is in fact true for all compact $L\subset \Omega_R\setminus K$ and this is an immediate corollary of a  Nguyen Thanh Van's result (\cite{nguyen}, page 228, see also \cite{nguyenzeriahi}, \cite{zeriahi} for ``pluricomplex versions''). At this point, it's not difficult to deduce 
$$\forall\,\varepsilon>0,\ \exists\, C_\varepsilon>0\ :\ \Vert\varphi_n\Vert_K\leq C_\varepsilon R^{n\varepsilon},\ \forall\,n\in\mathbb N,$$
and
$$\exists\, C>0\ :\ \Vert\varphi_n\Vert_L\geq C \cdot R^{n\frac{\alpha_L}{2}},\ \forall\,n\in\mathbb N.$$
It    remains to choose  $\varepsilon>0$ small enough so that $R^\varepsilon< R^{\frac{\alpha_L}{2}}$. Such a choice assures 
$$0\leq \lim_{n\to+\infty}\dfrac{\Vert \varphi_n \Vert_{K}}{ \Vert \varphi_n \Vert_{D(0,\rho_0)}}\leq \lim_{n\to+\infty} \left( R^{\varepsilon-\frac{\alpha_L}{2}}\right)^n=0.$$
So the sequence $\left(\frac{\Vert \varphi_n \Vert_{K}}{ \Vert \varphi_n \Vert_{D(0,\rho_0)}}\right)_n$ is bounded : by (8) there exists  $C>0$ such that 
$$\sum_{n\geq 1}
\vert a_n\vert \cdot \Vert\varphi_n\Vert_{K}\leq C\sum_{n\geq 1}\vert
a_n\vert\cdot\Vert\varphi_n\Vert_{D(0,\rho_0)},$$
which give us with (6) the estimation : $$\sum_{n\geq
1} \vert a_n\vert \cdot \Vert\varphi_n\Vert_{K}\leq 
C\big(\Vert f \Vert_{D(0,\rho_1)}-\vert f(0)\vert\big).$$ 
Finally, with the invariant Schwarz's lemma $(7)$ 
$$\sum_{n\geq
1} \vert a_n\vert \cdot \Vert\varphi_n\Vert_{K}\leq C\left( {\rho_1 \rho^{-1}+\vert f(0)\vert\over 1+\vert f(0)\vert
\rho_1 \rho^{-1}}-\vert f(0)\vert\right)= C\rho_1 \rho^{-1}\left({1-\vert f(0)\vert^2\over 1+\vert f(0)\vert
\rho_1 \rho^{-1}} \right)$$  which lead us to the main estimation
$$\sum_{n\geq
1} \vert a_n\vert \cdot \Vert\varphi_n\Vert_{K}\leq 2C\rho_1 \rho^{-1}(1-\vert f(0)\vert).$$
To conclude, let us choose $\rho$ large enough so that $2C\rho_1 \rho^{-1}\leq 1$, therefore, for any $R>1$ such that $D(0,\rho)\subset \Omega_R$ and 
$f=f(0)+\sum_{n\geq
1}a_n\varphi_n\in \mathscr O(\Omega_R)$, $f(\Omega_R)\subset \mathbb D$,  we have :
$$\vert f(0)\vert + \sum_{n\geq
1} \vert a_n\vert\cdot\Vert\varphi_n\Vert_{K}\leq 1.$$
Q.E.D.\hfill$\blacksquare$

\bigskip
Of course we must now   come back to the basis $(F_{K,n})_n$   :

\bigskip
\begin{lem} Let $\widetilde K\subset K$ be another compact, $(\varepsilon_n)_{n\geq 1}$ a complex sequence and suppose that there exists a constant $0<C<1$ such that
$$\begin{cases} \sup_{z\in \widetilde K}\vert \varphi_n(z)-\varepsilon_n\vert \leq C\cdot \Vert\varphi_n\Vert_K,\quad\forall\,n\in\mathbb N,&\qquad (9) \\
                \vert\varepsilon_n\vert\leq (1-C)\cdot \Vert\varphi_n\Vert_K,\quad\forall\,n\in\mathbb N.&\qquad (10)
                \end{cases}$$
Then the family  $( \widetilde K,\Omega, (\widetilde\varphi_n)_{n\geq 0})$ satisfies  Bohr's property with  $\widetilde\varphi_0\equiv 1,\ \widetilde\varphi_n:=\varphi_n-\varepsilon_n$.
\end{lem}

\bigskip
\noindent\textbf{Proof : } Let $f=a_0+\sum_{n\geq 1}a_n\varphi_n=a_0+\sum_{n\geq 1}a_n\epsilon_n+\sum_{n\geq 1}a_n(\varphi_n-\varepsilon_n)\in\mathscr O(\Omega)$ 
%\footnote{Il est bien de remarquer ‡ ce niveau que la condition (2) du lemme assure la convergence des deux sÈries $\sum_{n\geq 1}a_n\epsilon_n+\sum_{n\geq 1}a_n(\varphi_n-\varepsilon_n)$ il en rÈsulte alors facilement que la famille $(\widetilde\varphi_n)_{n\geq 0}$ est bien une base de $\mathscr O(\Omega)$.} 
and suppose that $\vert f\vert\leq 1$ on $\Omega$. We have to prove that 
$$\left\vert a_0+\sum_{n\geq 1}a_n\epsilon_n\right\vert+\sum_{n\geq 1}\vert a_n\vert\cdot\Vert \varphi_n-\varepsilon_n\Vert_{\widetilde K}\leq 1.$$
But :
$$\begin{aligned}
\left\vert a_0+\sum_{n\geq 1}a_n\epsilon_n\right\vert&+\sum_{n\geq 1}\vert a_n\vert\cdot\Vert \varphi_n-\varepsilon_n\Vert_{\widetilde K}
\leq \\
&\leq \left\vert a_0+\sum_{n\geq 1}a_n\epsilon_n\right\vert + C\cdot \sum_{n\geq 1}\vert a_n\vert\cdot\Vert \varphi_n\Vert_K \\
&\leq  \vert a_0\vert +\sum_{n\geq 1}\vert a_n\vert\cdot\vert \epsilon_n\vert+ C\cdot \sum_{n\geq 1}\vert a_n\vert\cdot\Vert \varphi_n\Vert_K\\
&\leq  \vert a_0\vert +(1-C)\sum_{n\geq 1}\vert a_n\vert\cdot\Vert \varphi_n\Vert_K+ C\cdot \sum_{n\geq 1}\vert a_n\vert\cdot\Vert \varphi_n\Vert_K\\
&\leq \vert a_0\vert +\sum_{n\geq 1}\vert a_n\vert\cdot\Vert \varphi_n\Vert_K\leq 1
\end{aligned}$$
Q.E.D.\hfill$\blacksquare$

\bigskip
\noindent\textbf{Proof of Theorem 3.1 : } Now let $K$ be a continuum, $\Omega_R,\ (R>1)$, a level set of the Green function of $K$ and fix $\widetilde K=\overline{\Omega_R }$. If $a\in\partial\Omega_R$ there exists (this is theorem 3.2) $R'>R$ such that the family $\left(\overline{\Omega_R }, \Omega_{R'}, (1, F_{\widetilde K, n}-F_{\widetilde K, n}(a))_{n\geq 0}\right)$ satisfies Bohr's property. Then for any function 
$$f=a_0+\sum_{n\geq 1} a_n ( F_{\widetilde K, n}-F_{\widetilde K, n}(a)) \in \mathscr O(\Omega_{R'}),$$  
such that $\vert f\vert \leq 1$ on $\Omega_{R'}$, we have  $$\vert a_0\vert+ \sum_{n\geq 1} \vert a_n\vert\cdot \Vert F_{\widetilde K, n}-F_{\widetilde K, n}(a)\Vert_{\overline{\Omega_R }}\leq 1.$$
\noindent But (\cite{suetin}, page 35) : 
$F_{\widetilde K, n}(z)=R^{-n}F_{K, n}(z)$ so
$$\begin{aligned} f(z)&=a_0+\sum_{n\geq 1} a_n \left( F_{\widetilde K, n}(z)-F_{\widetilde K, n}(a)\right)\\
 &= f(z)=a_0+\sum_{n\geq 1} a_n R^{-n}\left( F_{ K, n}(z)-F_{ K, n}(a)\right).
\end{aligned}$$
Because $R>1$, this immediately implies that the basis $(1, F_{K, n}-F_{K, n}(a))_{n\geq 0}$ satisfies Bohr's property on $(\overline{\Omega_R }, \Omega_{R'})$. If we  apply lemma 3.3 with $\varphi_n=F_{K, n}-F_{K, n}(a),\ a\in\partial \Omega_R$ and $-\varepsilon_n=F_{K, n}(a)$, the inequalities  (9) and (10) are :
$$\begin{aligned}
&(9')\qquad& \sup_{z\in K}\vert F_{K, n}(z)\vert \leq C\cdot \sup_{z\in\overline{\Omega_R}}\vert F_{K, n}(z)-F_{K, n}(a)\vert, \\
&(10')\qquad& \vert F_{K, n}(a)\vert \leq (1-C)\cdot \sup_{z\in\overline{\Omega_R}}\vert F_{K, n}(z)-F_{K, n}(a)\vert.
\end{aligned}$$
(where $C\in]0,1[$ is a constant). For all $n\in\mathbb N$ choose $a_n\in\partial\Omega_R$ such that $\Phi(a_n)=\theta_n\Phi(a)$ where $\theta_n$ is an $n$-root of $-1$ (remember that $\Phi(\partial\Omega_R)=C(0,R)$). So 
$$F_{K, n}^n(a_n)=\Phi^n(a_n)-E_{K, n}(a_n)=-\Phi(a)^n-E_{K, n}(a_n)$$
 and 
 $$F_{K, n}(a_n)-F_{K, n}(a)=-2\Phi(a)^n-\left[ E_{K, n}(a)+E_{K, n}(a_n)\right ].$$ 
 But because of inequality (3) in   paragraph 3.1, noting $r=1+\varepsilon_0$  :
$$\left\vert  E_{K, n}(a)+E_{K, n}(a_n)\right\vert \leq 2(1+\varepsilon_0)^n\dfrac{\texttt{lg}(\partial\Omega_{1+\varepsilon_0})}{\texttt{dist}(\partial\Omega_{1+\varepsilon_0},\partial\Omega_{R})},$$
for all $n\in\mathbb N$ and $R>1+\varepsilon_0$. Consequently :
$$\begin{aligned}  \sup_{z\in\overline{\Omega_R}}\vert F_{K, n}(z)-F_{K, n}(a)\vert &\geq \vert F_{K, n}(a_n)-F_{K, n}(a)\vert \\
&\geq 2R^n\left[1-\left(\dfrac{1+\varepsilon_0}{R}\right)^n\cdot\dfrac{\texttt{lg} (\partial\Omega_{1+\varepsilon_0})}{\texttt{dist}(\partial\Omega_{1+\varepsilon_0},\partial\Omega_{R})} \right] 
\end{aligned}$$
for all $n\in\mathbb N$ and $R>1+\varepsilon_0$. So, as long as we choose  $R$ large enough, say $R>R_0$, we can suppose that
%\footnote{En effet, le terme entre crochets vÈrifie pour $R>R'>1+\varepsilon_0$ :  $$1>1-\left(\dfrac{1+\varepsilon_0}{R}\right)^n\cdot\dfrac{\texttt{lg} (\partial\Omega_{1+2\varepsilon_0})}{\texttt{dist}(\partial\Omega_{1+\varepsilon_0},\partial\Omega_{R})}>1-\left(\dfrac{1+\varepsilon_0}{R_1}\right)\cdot\dfrac{\texttt{lg} (\partial\Omega_{1+2\varepsilon_0})}{\texttt{dist}(\partial\Omega_{1+\varepsilon_0},\partial\Omega_{R_1})}\underset{R_1\to\infty}{\longrightarrow} 1.$$
%on peut donc  le minorer par $3/4$ pour tout $n\in\mathbb N$ et $R>R_1$}
$$  \sup_{z\in\overline{\Omega_R}}\vert F_{K, n}(z)-F_{K, n}(a)\vert \geq \dfrac{3}{2}R^n,\quad \forall\,n\in\mathbb N,\ R>R_0.\leqno{(11)}$$
Because of (4)  :
$$\vert F_{K,n}(a)\vert \leq R^n\left[ 1+\left(\dfrac{1+\varepsilon_0}{R}\right)^n\cdot\dfrac{\texttt{lg} (\partial\Omega_{1+\varepsilon_0})}{\texttt{dist}(\partial\Omega_{1+\varepsilon_0},\partial\Omega_{R})} \right]$$
for all $n\in\mathbb N$, $R>1+\varepsilon_0$. Because the term in between the brackets satisfies :  
$$\begin{aligned}1&\leq 1+\left(\dfrac{1+\varepsilon_0}{R}\right)^n\cdot\dfrac{\texttt{lg} (\partial\Omega_{1+\varepsilon_0})}{\texttt{dist}(\partial\Omega_{1+\varepsilon_0},\partial\Omega_{R})}\\
&\leq 1+\left(\dfrac{1+\varepsilon_0}{R_1}\right)\cdot\dfrac{\texttt{lg} (\partial\Omega_{1+\varepsilon_0})}{\texttt{dist}(\partial\Omega_{1+\varepsilon_0},\partial\Omega_{R_1})}\underset{R_1\to\infty}{\longrightarrow} 1
\end{aligned}$$
 for all $R>R_1>1+\varepsilon_0$ ; it is less than $5/4$ for all $n\in\mathbb N$ and $R>R_1$ where $R_1$ is choosen large enough ; i.e.
$$\vert F_{K,n}(a)\vert \leq \dfrac{5}{4}\cdot R^n,\quad \forall\,n\in\mathbb N,\ R>R_1.$$
It follows from (11) that
$$\vert F_{K,n}(a)\vert \leq \dfrac{5}{6}\cdot \dfrac{3}{2}\cdot R^n \leq \dfrac{5}{6}\sup_{z\in\overline{\Omega_R}}\vert F_{K, n}(z)-F_{K, n}(a)\vert,\ \forall\,n\in\mathbb N,\ R>R_2:=\max\{ R_0, R_1\}.$$
So we have proved  inequality  (10') with $C=1/6$. 
 Finaly, still because of (4) :
$$\begin{aligned}
\qquad \sup_{z\in K}\vert F_{K, n}(z)\vert &\leq \sup_{z\in\overline{\Omega_{1+2\varepsilon_0}}}\vert F_{K, n}(z)\vert \\
&\leq (1+2\varepsilon_0)^n\cdot \left[ 1+\left(\dfrac{1+\varepsilon_0}{1+\varepsilon_0}\right)^n\cdot\dfrac{\texttt{lg} (\partial\Omega_{1+\varepsilon_0})}{\texttt{dist}(\partial\Omega_{1+2\varepsilon_0},\partial\Omega_{1+\varepsilon_0})} \right]\\
&\leq A (1+2\varepsilon_0)^n,\quad\forall\,n\in\mathbb N
\end{aligned}$$ 
where $A$ is a constant strictly larger than $1$. Given  $A>1$ being fixed it is easy to deduce that for any $R>R_3$ :
$$\sup_{z\in K}\vert F_{K, n}(z)\vert \leq A (1+2\varepsilon_0)^n \leq \dfrac{R^n}{4},\quad \forall\, n\in\mathbb N.$$
So because of (11)
$$\sup_{z\in K}\vert F_{K, n}(z)\vert \leq \dfrac{1}{6}\cdot\dfrac{3}{2}R^n \leq \dfrac{1}{6} \sup_{z\in\overline{\Omega_R}}\vert F_{K, n}(z)-F_{K, n}(a)\vert$$
for all $n\in\mathbb N$ et $R > \max\{ R_3, R_2\}$. This is formula (9') with $C=1/6$, so we can apply lemma 3.3 and deduce that the family $(K,\Omega_R, (F_{K, n})_{n\geq 0})$ satisfies  Bohr's phenomenon for all $R$ large enough :  theorem 3.1. is proved.\hfill$\blacksquare$

\end{document}